\documentclass[nonblindrev]{informs3_noinforms}

\DoubleSpacedXI 

\usepackage{endnotes}
\let\footnote=\endnote

\usepackage{enumerate, float, tikz}

\newcommand{\E}{\bbe}

\newcommand{\p}{\mathbb{P}}

\def\bbr{{\Bbb{R}}} 
\def\bbe{{\Bbb{E}}} 

\def\bbp{{\Bbb{P}}}

\renewcommand{\theequation}{\thesection.\arabic{equation}}

\newtheorem{observation}{Observation}

\usepackage{natbib}
 \bibpunct[, ]{(}{)}{,}{a}{}{,}%
 \def\bibfont{\small}%
\TheoremsNumberedThrough     
\ECRepeatTheorems

\EquationsNumberedThrough    

\begin{document}



\RUNTITLE{Optimality gap of constant-order policies decays exponentially in the lead time for lost sales models}
\TITLE{Optimality gap of constant-order policies decays exponentially in the lead time for lost sales models}

\ARTICLEAUTHORS{%
\AUTHOR{\bf Linwei Xin}
\AFF{Department of Industrial and Enterprise Systems Engineering,
University of Illinois at Urbana-Champaign,
Urbana, IL 61801}
\AUTHOR{
\bf David A. Goldberg}
\AFF{School of Industrial \& Systems Engineering,
Georgia Institute of Technology,
Atlanta, GA 30332}

\RUNAUTHOR{{Xin and Goldberg}}

} 

\ABSTRACT{
\indent Inventory models with lost sales and large lead times have traditionally been considered intractable due to the curse of dimensionality. Recently, Goldberg and co-authors laid the foundations for a new approach to solving these models, by proving that as the lead time grows large, a simple constant-order policy is asymptotically optimal. However, the bounds proven there require the lead time to be very large before the constant-order policy becomes effective, in contrast to the good numerical performance demonstrated by Zipkin even for small lead time values.  In this work, we prove that for the infinite-horizon variant of the same lost sales problem, the optimality gap of the same constant-order policy actually converges \emph{exponentially fast} to zero, with the optimality gap decaying to zero at least as fast as the exponential rate of convergence of the expected waiting time in a related single-server queue to its steady-state value. We also derive simple and explicit bounds for the optimality gap, and demonstrate good numerical performance across a wide range of parameter values for the special case of exponentially distributed demand.  Our main proof technique combines convexity arguments with ideas from queueing theory.
}



\maketitle

%


\section{Introduction}
It is well-known that there is a fundamental dichotomy in the theory of inventory models, depending on the fate of unmet demand.
If unmet demand remains in the system and can be met at a later time, we say the system exhibits \emph{backlogged demand}; if unmet demand is lost to the system, we say the system exhibits \emph{lost sales}.  Which of these assumptions is appropriate depends heavily on the application of interest.  For example, in many retail applications demand can be met by a competing supplier, making lost sales a more appropriate assumption (cf. \cite{BV}).

A second important feature of many inventory models, intimately related to the above dichotomy, is that of positive lead times, i.e. settings in which there is a multi-period delay between when an order for more inventory is placed and when that order is received.  This feature leads to an enlarged state-space (growing linearly with the lead time), to track all orders already placed but not yet received, i.e. the pipeline vector.  It is a classical result that models with backlogged demand remain tractable even in the presence of positive lead times.  Namely, it can be proven that a so-called base-stock (i.e. order-up-to) policy, based only on the total inventory position (i.e. sum of the current inventory and all orders in the pipeline vector), is optimal in this setting (cf. \cite{Scarf60,Iglehart,Veinott66}).  However, it is known that such simple policies are no longer optimal for models with lost sales and positive lead times (cf. \cite{KS}).  \emph{For over fifty years, inventory models with lost sales and positive lead times were generally considered intractable, as the primary solution method (dynamic programming) suffered from the curse of dimensionality as the lead time grew.}  As noted in \cite{BV}, this has led to many researchers using models with backlogging as approximations for settings in which a lost sales assumption is more appropriate, which may lead to very suboptimal policies.

Although the optimal policy for lost-sales models with positive lead times remains poorly understood, the model has been studied now for over fifty years, and we refer to \cite{BV} and \cite{GKLSS} for a comprehensive review.  Due to the difficulty of computing the optimal policy, there has been considerable focus on understanding structural properties of an optimal policy, and analyzing heuristics.  In particular, convexity results were obtained in \cite{KS}, \cite{Morton69}, and \cite{Zipkina}.  \cite{JSS07} compared the optimal costs between the backlogged and lost sales systems with identical problem parameters, and \cite{HJMR} further proved that the base-stock policy was asymptotically optimal as the lost-sales penalty became large compared to the holding cost, with similar results also derived in \cite{LuSqYa12}.
\cite{LJN} proposed the family of so called dual-balancing policies, and proved that the cost incurred by such a policy was always within a factor of $2$ of optimal.  Another recent line of research yields efficient approximation algorithms \emph{for any fixed lead time} by carefully approximating the associated dynamic programs (cf. \cite{Halman09, Halman12,Chen12a}).  \emph{Despite this progress, the aforementioned work leaves open the problem of deriving efficient algorithms with arbitrarily small error, when the lead time is large.}

A very simple and natural policy, which will be the subject of our own investigations, is the so-called \emph{constant-order policy}, which places the same order in every period, independent of the state of the system.  Perhaps surprisingly, \cite{Reiman} proved that for lost sales inventory models with positive lead times, sometimes the best constant-order policy outperforms the more sophisticated base-stock policy, and performed a detailed analysis under a certain asymptotic scaling.  This phenomena was further illuminated by the computational study of \cite{Zipkinb}, which confirmed that in several scenarios the constant-order policy performed favorably.  These observations were recently given a solid theoretical foundation by \cite{GKLSS}, who proved that for lost sales inventory models with positive lead times, as the lead time grows with all other parameters remaining fixed, the best constant-order policy is in fact \emph{asymptotically optimal}.  This is quite surprising, as the policy is so simple, and performs nearly optimally exactly in the setting which had stumped researchers for over fifty years.   However, the bounds proven there are impractical, requiring the lead time to be very large before the constant-order policy becomes nearly optimal, e.g. requiring a lead time which is $\Omega(\epsilon^{-2})$ to ensure a $(1 + \epsilon)$-approximation guarantee, and involving a massive prefactor.  The authors note that the numerical experiments of \cite{Zipkinb} suggest that the constant-order policy performs quite well even for small lead times, and pose closing this gap (thus making the results practical) as an open problem.  The authors also point out that if one could prove that the constant-order policy performs well for small to moderate lead times, this would open the door for the creation of practical hybrid algorithms, which solve large dynamic programs when the lead time is small, and transition to simpler algorithms for larger lead times.

\subsection{Our contributions}
In this work, we make significant progress towards solving this open problem and closing this gap.  For the infinite-horizon variant of the finite-horizon problem considered by \cite{GKLSS}, we prove that the optimality gap of the same constant-order policy actually converges \emph{exponentially fast} to zero, i.e. we prove that a lead time which is $O\big(\log(\epsilon^{-1})\big)$ suffices to ensure a $(1 + \epsilon)$-approximation guarantee.  We demonstrate that the corresponding rate of exponential decay is at least as fast as the exponential rate of convergence of the expected waiting time in a related single-server queue to its steady-state value, which we prove to be monotone in the ratio of the lost-sales penalty to the holding cost.  We also derive simple and explicit bounds for the optimality gap.  For the special case of exponentially distributed demand, we further compute all expressions appearing in our bound in closed form, and numerically evaluate them, demonstrating good performance for a wide range of parameter values.  Our main proof technique combines convexity arguments with ideas from queueing theory, and is simpler than the coupling argument of \cite{GKLSS}.

\subsection{Outline of paper}
The rest of the paper is organized as follows. We formally define the inventory problem of interest in Section \ref{ls-sec-prob}, and describe the constant-order policy in Section \ref{ls-sec-constant}. We state our main results in Section \ref{ls-subsec-main}, and prove our main results in Section\ \ref{ls-sec-proof}.  We summarize our results and propose directions for future research in Section\ \ref{ls-sec-conc}.  A technical appendix is provided in Section\ \ref{ls-sec-appen}, which contains several proofs from throughout the paper, as well as an in-depth analysis of our bounds (both analytic and numeric) for the illustrative special case of exponentially distributed demand.  

\section{Main results}\label{ls-sec-main}
\subsection{Model description and problem statement}\label{ls-sec-prob}
In this section, we formally define our lost-sales inventory optimization problem.  Let $\lbrace D_t, t \geq 1 \rbrace$ be a sequence of independent and identically distributed (i.i.d.) demand realizations, distributed as the non-negative random variable (r.v.) $D$ with distribution ${\mathcal D}$, which we assume to have finite mean, and (to rule out certain trivial degenerate cases) to have strictly positive variance.  Let $L \geq 1$ be the deterministic lead time, and $h,p$ ($>0$) be the unit holding cost and lost-sales penalty respectively. In addition, let $I_t$ denote the on-hand inventory, and $\mathbf{x}_t = (x_{1,t},\ldots, x_{L,t})$ denote the pipeline vector of orders placed but not yet delivered, at the beginning of time period $t$, where $x_{i,t}$ is the order to be received in period $i + t - 1$. The ordered sequence of events in period $t$ is then as follows.
\begin{itemize}
\item A new amount of inventory $x_{1,t}$ is delivered and added to the on-hand inventory.
\item A new order is placed.
\item The demand $D_t$ is realized.
\item Costs for period $t$ are incurred, and the on-hand inventory and pipeline vector are updated.
\end{itemize}
Note that the on-hand inventory is updated according to $I_{t+1} =\max(0, I_t + x_{1,t} - D_t)$, and the pipeline vector is updated such that (s.t.) $x_{1,t}$ is removed, $x_{i,t+1}$ is set equal to $x_{i+1,t}$ for $i \in [1,L-1]$, and $x_{L,t+1}$ is set equal to the new order placed.  We now formalize the family of admissible policies $\Pi$, which will determine the new order placed, mirroring the discussion in \cite{HJN}.  An admissible policy $\pi$ consists of a sequence of measurable maps $\lbrace f^{\pi}_t, t \geq 1 \rbrace$, where each $f^{\pi}_t$ is a deterministic measurable function with domain ${\mathcal R}^{+, L + 1}$ and range ${\mathcal R}^+$.  In that case, for a given policy $\pi$, the order placed in period $t$ equals $f^{\pi}_t(\mathbf{x}_t,I_t)$; and $\Pi$ denotes the family of all such admissible policies $\pi$.  We let $C_t \stackrel{\Delta}{=}\ h\left(I_t + x_{1,t} - D_t\right)^+  + p \left(I_t + x_{1,t} - D_t\right)^-$ denote the sum of the holding cost and lost-sales penalty in time period $t$,  
where $x^+\stackrel{\Delta}{=} \max(x, 0)$, $x^-\stackrel{\Delta}{=} \max(-x, 0)$.   Let $\mathbf{0} (\mathbf{1})$ denote the all-zeros (all-ones) vector, where the dimension is to be inferred from context.  As we will be primarily interested in minimizing long-run average costs, and for simplicity, we suppose that the initial conditions are such that the initial inventory is 0, and the initial pipeline vector is empty, i.e. $I_1 = 0, \mathbf{x}_1 = \mathbf{0}$.  For a given policy $\pi \in \Pi$, we let $C^{\pi}_t$ denote the corresponding cost incurred in period $t$, and $C(\pi)\stackrel{\Delta}{=}\  \limsup_{T\to \infty}\frac{\sum_{t=1}^{T} \bbe\left[C^{\pi}_t \right]}{T}$ denote the long-run average cost incurred by $\pi$.  Then the corresponding infinite-horizon lost-sales inventory optimization problem which we will consider (identical to the problem considered in \cite{HJN}) is
\begin{equation} \label{prob1}
\textrm{OPT}(L) \stackrel{\Delta}{=}\  \inf_{\pi\in \Pi} C(\pi).
\end{equation}
For a given policy $\pi \in \Pi$, let $\mathbf{x}^{\pi}_t (I^{\pi}_t)$ denote a random vector (variable) distributed as the pipeline vector (inventory) at the start of period $t$ under policy $\pi$.  Further suppose that $\lbrace (\mathbf{x}^{\pi}_t, I^{\pi}_t), t \geq 1 \rbrace$ are all constructed on a common probability space, and have the appropriate joint distribution induced by the operation of $\pi$ over time.  Recall that a policy $\pi \in \Pi$ is said to be \emph{stationary} if $f^{\pi}_t$ is independent of $t$, i.e. $f^{\pi}_t = f^{\pi}_1$ for all $t$, in which case we denote this common function by $f^{\pi}$.  In that case, it follows immediately from the results of \cite{HJN} and \cite{Zipkina} that for all $L \geq 1$, there exists a stationary optimal policy $\pi^{*,L}$ for Problem (\ref{prob1}) whose induced inventory level and pipeline vector is uniformly bounded over time.  More formally, for all $L \geq 1$, there exists a stationary policy $\pi^{*,L} \in \argmin_{\pi \in \Pi} C(\pi)$ and finite positive constant $M_L$ s.t. with probability 1 (w.p.1), for all $t \geq 1$,
\begin{equation}\label{isbounded}
\sum_{i=1}^L x^{\pi^{*,L}}_{i,t} + I^{\pi^{*,L}}_t < M_L.
\end{equation}
As a notational convenience, we denote $\mathbf{x}^{\pi^{*,L}}_t$ by $\mathbf{x}^{*,L}_t$, and $I^{\pi^{*,L}}_t$ by $I^{*,L}_t$.

\subsection{Constant-order policy}\label{ls-sec-constant}
In this section, we formally define the constant-order policy, and characterize the best constant-order policy.  As a notational convenience, let us define all empty sums to equal zero, let $e$ denote Euler's number, $\log(x)$ denote the natural logarithm of $x$, $\frac{1}{\infty}$ denote $0$, $\frac{1}{0}$ denote $\infty$, $\log(\infty)$ denote $\infty$, and $\mathbb{I}(A)$ denote the indicator of the event $A$.  For any $r \in [0, \E[D])$, the constant-order policy $\pi_r$ is the policy that places the constant order $r$ in every period.  It is well-known (cf. \cite{Asmussen}) that for any such $r \in [0,\E[D])$, $\lbrace I^{\pi_r}_t, t \geq 1 \rbrace$ converges in distribution, and in expectation, to a limiting r.v. with finite mean, which we denote by $I^r_\infty$.  $I^r_{\infty}$ has the same distribution as the steady-state waiting time in the corresponding $GI/GI/1$ queue with interarrival distribution ${\mathcal D}$ and processing time distribution the constant $r$.  For two r.v.s $X,Y$, let $X \sim Y$ denote equivalence in distribution between $X$ and $Y$.  In that case, it is well-known (cf. \cite{Asmussen}) that $I^r_\infty \sim \sup_{j \geq 0} \left(j r - \sum_{i=1}^j D_{i}\right)$, and that for any $r \in [0,\E[D])$ 
\begin{equation}\label{lseqconst}
C(\pi_r) = \ h \E\left[\sup_{j \geq 0} \left(j r - \sum_{i=1}^j D_{i}\right)\right] + p \E[D] - p r,
\end{equation}
independent of the lead time $L$.  As it is easily verified from (\ref{lseqconst}) that $C(\pi_r)$ is a convex function of $r$ on $[0,\E[D])$, to find the best possible constant-order policy, it suffices to select the $r$ minimizing this one-dimensional convex function over the compact set $[0,\E[D]]$.  We note that the existence of at least one such optimal $r$ follows from the well-known properties of convex optimization over a compact set, and that the set of all such optimal solutions must be bounded away from $\E[D]$, since by assumption $h > 0$, and $\lim_{r \uparrow \E[D]} \bbe[I^r_\infty] = \infty$ (since $D$ has strictly positive variance, cf. \cite{Asmussen}).  Let $r_{\infty} \in\ \arg\min_{0\leq r\leq \bbe[D]} C(\pi_r)$ denote the infimum of this set of optimal constant-order levels, in which case the best constant-order policy will refer to $\pi_{r_{\infty}}$.

\subsection{Main results}\label{ls-subsec-main}
In this section, we present our main results, demonstrating that the optimality gap of the best constant-order policy decays exponentially in the lead time.  For $\theta \geq 0$, let us define 
$$\phi(\theta) \stackrel{\Delta}{=} \exp(\theta r_{\infty}) \bbe[\exp(- \theta D)]\ \ \ ,\ \ \ \gamma \stackrel{\Delta}{=} \inf_{\theta \geq 0} \phi(\theta),$$
and $\vartheta \in \argmin_{\theta \geq 0} \phi(\theta)$ denote the supremum of the set of minimizers of $\phi(\theta)$, 
where we define $\vartheta$ to equal $\infty$ if the above infimum is not actually attained.  Note that $\phi(\theta)$ is a continuous and convex function of $\theta$ on $(0,\infty)$, and right-continuous function of $\theta$ at $0$.    
In addition, it follows from \cite{Folland} Theorem 2.27 that $\phi(\theta)$ is right-differentiable at zero, with derivative equal to $r_{\infty} - \bbe[D]$.  As $r_{\infty} < \bbe[D]$, we conclude from the definition of derivative and a straightforward contradiction argument that $\vartheta > 0$ (i.e. $\vartheta$ is strictly positive), and $\gamma \in [0,1)$ (i.e. $\gamma$ is strictly less than 1).  It follows from the celebrated Cram\'{e}r's Theorem, and more generally the theory of large deviations, that up to exponential order $\p( k r_{\infty} \geq \sum_{i=1}^k D_i )$ decays like $\gamma^k$ as $k \rightarrow \infty$ (cf. \cite{DS}).  Furthermore, as we will explore in detail later in the proof of our main result, $\gamma$ corresponds (again up to exponential order) to the rate at which the expected waiting time in an initially empty single-server queue, with inter-arrival distribution $D$ and processing time distribution (the constant) $r_{\infty}$, converges to its steady-state value (cf. \cite{Kingman}).  
Let $g \stackrel{\Delta}{=} \inf_{x \in {\mathcal R}} \E\big[ h(x-D)^+ + p(x - D)^-\big] > 0$, 
and $Q$ denote the $\frac{p}{p + h}$ quantile of the demand distribution, where we note that $Q$ is the optimal order quantity in the single-stage newsvendor problem and $\textrm{OPT}(L) \geq g$ (cf. \cite{Zipkinc}).  Then our main result is as follows.
\begin{theorem}[Exponential convergence of constant-order policy to optimality]\label{lsthm-main}
For all $L \geq 1$, 
\begin{equation}\label{ourmainbound1}
\frac{C(\pi_{r_{\infty}})}{\textrm{OPT}(L)} \leq 1 + h \big( (1 - \gamma) g \big)^{-1} \bigg( \E[D] - r_{\infty} + \big(e \vartheta (L+1) \big)^{-1} \bigg) \gamma^{L+1}.
\end{equation}
\end{theorem}

Our results prove that for the corresponding infinite-horizon problem, the optimality gap of the constant-order policy converges exponentially fast to zero.  In particular, a lead time which is $O\big(\log(\epsilon^{-1})\big)$ suffices to ensure a $(1 + \epsilon)$-approximation guarantee.  This contrasts with the bounds of \cite{GKLSS}, which had an inverse polynomial dependence on $\epsilon$.  Furthermore, our explicit bounds are much tighter than those of \cite{GKLSS}.  This takes a large step towards answering several open questions posed in \cite{GKLSS} with regards to deriving bounds tight enough to be useful in practice.  We note that since our results only hold for the infinite-horizon problem, and will use critically certain relationships which only hold in this regime, our results are not directly comparable to those of \cite{GKLSS}, whose bounds also hold for finite-horizon problems.  Closing this gap, and proving tighter bounds for the finite-horizon problem, remains an interesting open question.  To further illustrate the performance of our bounds, in the technical appendix Section\ \ref{ls-sec-appen}
we explicitly compute (in closed form) the right-hand side of (\ref{ourmainbound1}) for the special case that demand is exponentially distributed, and demonstrate good numerical performance for a wide range of parameter values.

\subsubsection{Impact of the ratio $\frac{p}{h}$ on our bounds}\label{ls-sec-impact}
In this section, we discuss the dependence of our demonstrated exponential rate of convergence $\gamma$ on the ratio of the lost-sales penalty to the holding cost.  In particular, we show that $\gamma$ is non-decreasing in $\frac{p}{h}$.  Note that by a simple scaling argument, for any fixed demand distribution ${\mathcal D}$, $r_{\infty}$ (and thus $\gamma$) is a function of $\frac{p}{h}$ only, as opposed to the particular values of $p,h$.  To make this dependence explicit, let $\gamma(\varrho)$ denote the value of $\gamma$ when $\frac{p}{h} = \varrho$ (with the dependence on ${\mathcal D}$ implicit).
\begin{lemma}\label{lsprop-main1}
For any fixed demand distribution ${\mathcal D}$, $\gamma(\varrho)$ is non-decreasing in $\varrho$.
\end{lemma}
We include a proof of Lemma\ \ref{lsprop-main1} in the technical appendix Section\ \ref{ls-sec-appen}.  This result suggests that the optimality gap of the constant-order policy may be larger when $\frac{p}{h}$ is large.  Interestingly, this is exactly the regime in which \cite{HJMR} proved that order-up-to policies are nearly optimal.  More formally understanding this connection remains an interesting open question.

\section{Proof of Theorem \ref{lsthm-main}}\label{ls-sec-proof}
In this section, we complete the proof of our main result Theorem \ref{lsthm-main}.  Let $\delta_{i,j}$ equal $1$ if $i=j$ and 0 otherwise.  Our proof begins by relating the ``long-run behavior" of $\pi^{*,L}$ to a certain constant-order policy.  At a high level, we will combine Lindley's recursion and convexity with the fact that the r.v.s corresponding to (appropriately defined stationary versions of) the different components of the pipeline vector (under $\pi^{*,L}$) have the same mean, which will (approximately) coincide with the constant in our constant-order policy.  However, this program immediately encounters a technical problem.  Namely, the natural way to analyze the ``long-run behavior" of $\pi^{*,L}$ is through the steady-state distribution of the Markov chain induced by the policy $\pi^{*,L}$, i.e. $\lbrace (\mathbf{x}^{*,L}_t, I^{*,L}_t), t \geq 1 \rbrace$.  However, it is not obvious that this steady-state exists, i.e. that $\lbrace (\mathbf{x}^{*,L}_t, I^{*,L}_t), t \geq 1 \rbrace$ converges weakly to a limiting stationary vector as $t \rightarrow \infty$.  To overcome this, we first observe that we will not actually need a random vector which is truly the steady-state of the aforementioned Markov chain (which in principle may not exist), but only need to demonstrate the existence of a random vector which has several properties that we would want such a steady-state (if it existed) to have.  We now show the existence of such a random vector.  We note that although closely related questions have been studied in the Markov decision process (MDP) literature (cf. \cite{Ara13}), and perturbative approaches similar to the approach we take in our own proof are in general well-known (cf. \cite{Filar}), to the best of our knowledge the desired result does not follow directly from any results appearing in the literature.  As such, we include a proof for completeness in the technical appendix Section\ \ref{ls-sec-appen}.  

\begin{theorem}\label{proxychain1}
One may construct  an $L$-dimensional random vector $\mathbf{\chi}^*$, r.v. $\mathcal{I}^*$, and $\lbrace D_i, i \geq 1 \rbrace$ on a common probability space s.t. the following are true.
\begin{enumerate}[(i)]
\item \label{enum1} $(\mathbf{\chi}^*, \mathcal{I}^*)$ is w.p.1 non-negative, has finite mean, and is independent of $\lbrace D_i, i \geq 1 \rbrace$.
\item \label{lseq-30}
$\E[\mathcal{I}^*] =\ \bbe\left[\max_{j=0,\ldots,L}\left(\sum_{i=1}^j(\mathbf{\chi}^*_{L+1-i}- D_{L+1-i}) + \delta_{j,L} \mathcal{I}^* \right) \right].$
\item \label{enum3} $\E[\mathbf{\chi}^*_i] = \E[\mathbf{\chi}^*_1] \leq \E[D]$ for all $i$.
\item \label{enum4} $h \E[{\mathcal I}^*] +  p \bbe[D_1]- p \bbe[\mathbf{\chi}^*_1] = \textrm{OPT}(L).$
\end{enumerate}
\end{theorem}
We note that in the above, $\mathbf{\chi}^*_1$ corresponds intuitively to (a stationary version) of the first component of the pipeline vector under $\pi^{*,L}$, i.e. the amount of inventory received in any single period.  $\mathbf{\chi}^*_1$ \emph{does not} correspond to any notion of inventory position (i.e. some aggregate measure of the total inventory in the pipeline), nor does $\mathbf{\chi}^*_i$ for any $i$.  Indeed, although certain previous papers on lost sales models sometimes perform their analysis on a transformed set of variables which do correspond to inventory positions (cf. \cite{Zipkina}), we \emph{do not} use any such transformation here.

Next, we observe that the right-hand side of Theorem\ \ref{proxychain1}.(\ref{lseq-30}) is a jointly convex function of $\chi^*$ and $\mathcal{I}^*$, which will allow us to apply the multi-variate Jensen's inequality (cf. \cite{Dudley}).  More formally, for fixed $\mathbf{d} \in \bbr^L$ and $\chi_1,\ldots,\chi_L,I \in {\mathcal R}$, let us define $\mathfrak{f}_{\mathbf{d}}\left(\chi_1,\ldots,\chi_L,I\right) \stackrel{\Delta}{=} \max_{j=0,\ldots,L}\left(\sum_{i=1}^j(\chi_{L+1-i}- d_{L+1-i}) +\delta_{j,L} I\right)$.

\begin{observation}\label{lslemma-jensen}
For each fixed $\mathbf{d} \in \bbr^L$, $\mathfrak{f}_{\mathbf{d}}\left(\chi_1,\ldots,\chi_L,I\right)$ is a jointly convex function of  $\left(\chi_1,\ldots,\chi_L,I\right)$ over $\bbr^{L+1}$.  Combining with Theorem\ \ref{proxychain1}.(\ref{enum1}) - (\ref{enum3}), the multi-variate Jensen's inequality, and the i.i.d. property of $\lbrace D_i, i \geq 1 \rbrace$, we conclude that
$$\E[\mathcal{I}^*] \geq \bbe\left[\max_{j=0,\ldots,L}\left(j \E[\chi^*_1] - \sum_{i=1}^j D_{i} +\delta_{j,L} \E[\mathcal{I}^*] \right)\right].$$
\end{observation}
We note that the observation of such a convexity in terms of the on-hand inventory and pipeline vector is not new.  Indeed, the so called $L$-natural-convexity of the relevant cost-to-functions has been studied extensively (cf. \cite{KS, Morton69, Zipkina, Chen12a}), and used to obtain both structural results and algorithms.  In contrast, here we use convexity to relate the expected inventory under an optimal policy to the expected inventory under a particular constant-order policy, intuitively that which orders $\bbe[\chi^*_1]$ in every period.  Very similar ideas and arguments have appeared previously in the queueing theory literature, to demonstrate the extremality (with regards to expected waiting times) of certain queueing systems with constant service (or inter-arrival) times (cf. \cite{Humblet, Hajek}).  

Before proceeding, let us define several additional notations.  In particular, for $r \in \big[0,\E[D]\big]$ and $L \geq 1$, let 
$$I^r_L\stackrel{\Delta}{=}\ \max_{j=0,\ldots,L}\left(jr - \sum_{i=1}^j D_{i}\right)\ \ \ ,\ \ \ C_L(r)
\stackrel{\Delta}{=}\ h\bbe[I^r_L]+ p \bbe[D]- p r,$$
and $r_L \in \arg\min_{0\leq r\leq \bbe[D]} C_L(r)$ denote the infimum of the set of minimizers of $C_L(r)$.  We note that $I^r_L$ is distributed as the waiting time of the $L$-th customer in the corresponding $GI/GI/1$ queue (initially empty) with interarrival distribution $D$ and processing time $r$.  We also note that for $r \in [0,\E[D])$, $I^r_{\infty}$ is the weak limit, as $L \rightarrow \infty$, of $I^r_L$.  Similarly, $C(\pi_r) = \lim_{L \rightarrow \infty} C_L(r)$, and $C_L(r)$ is monotone increasing in $L$.

We now combine Theorem\ \ref{proxychain1} with Observation\ \ref{lslemma-jensen} to bound the optimality gap of the constant-order policy.

\begin{lemma}\label{lslemma-0}
$\textrm{OPT}(L) \geq C_L(r_L)$, and
\begin{equation}\label{tobound2}
C(\pi_{r_{\infty}}) - \textrm{OPT}(L) \leq\ h\left(\bbe\left[I^{r_\infty}_\infty\right]- \bbe\left[I^{r_\infty}_L\right]\right)+   h\left(\bbe\left[I^{r_\infty}_L\right]- \bbe\left[I^{r_L}_L\right]\right) - p\left(r_\infty - r_L\right).
\end{equation}
\end{lemma}

\proof{Proof}
Combining Observation\ \ref{lslemma-jensen} with the nonnegativity of $\E[\mathcal{I}^*]$, we conclude that $\E[\mathcal{I}^*] \geq \bbe\left[I^{\E[\chi^*_1]}_L\right]$.  Thus by Theorem\ \ref{proxychain1}.(\ref{enum4}), $\textrm{OPT}(L) \geq h\bbe\left[I^{\E[\chi^*_1]}_L\right] + p \E[D] - p \E[\chi^*_1]$.  Combining with Theorem\ \ref{proxychain1}.(\ref{enum3}) and the definition of $r_L$, we conclude that $\textrm{OPT}(L) \geq C_L(r_L)$.  It then follows from (\ref{lseqconst}) that
\[
\begin{aligned}
   C\left(\pi_{r_\infty}\right) - \textrm{OPT}(L) \leq&\ \left(h\bbe\left[I^{r_\infty}_\infty\right]+ p \bbe[D]- p r_\infty\right)  -  \left(h \bbe\left[I^{r_L}_L\right]+ p \bbe[D]- p r_L\right)\\
   =&\ h\left(\bbe\left[I^{r_\infty}_\infty\right]- \bbe\left[I^{r_\infty}_L\right]\right)+   h \left(\bbe\left[I^{r_\infty}_L\right]- \bbe\left[I^{r_L}_L\right]\right) - p\left(r_\infty - r_L\right),
\end{aligned}
\]
completing the proof.  $\Halmos$
\endproof
We proceed by bounding the terms appearing in the right-hand side of (\ref{tobound2}) separately.  We begin by recalling a classical result of \cite{Kingman}, which uses the celebrated Spitzer's identity to bound the difference between the expected waiting time of the $L$th job to arrive to a single-server queue (initially empty), and the steady-state expected waiting time.  As this difference is exactly $\bbe\left[I^{r}_{\infty}\right]- \bbe\left[I^{r}_L\right]$, the result will allow us to bound the relevant term of (\ref{tobound2}).  We state Kingman's results as customized to our own setting, notations, and assumptions.
\begin{lemma}[Theorems 1 and 4 of \cite{Kingman}]\label{lslemma-9}
For all $r \in \big[ 0,\bbe[D] \big]$ and $L \geq 1$,
$$
   \bbe[I^r_L]=\ \sum_{n=1}^L \frac{1}{n} \bbe\left[\left(nr - \sum_{i=1}^n D_{i}\right)^+\right].
$$
If in addition $r < \bbe[D]$, then
$$
\bbe[I^r_\infty]=\ \sum_{n=1}^\infty \frac{1}{n} \bbe\left[\left(nr - \sum_{i=1}^n D_{i}\right)^+\right].
$$
Also,
\[
   \bbe\left[I^{r_{\infty}}_{\infty}\right]- \bbe\left[I^{r_{\infty}}_L\right]\leq\ \big( (1 - \gamma) e \vartheta (L+1) \big)^{-1} \gamma^{L+1}.
\]
\end{lemma}

To bound the remaining term $h\left(\bbe\left[I^{r_\infty}_L\right]- \bbe\left[I^{r_L}_L\right]\right) - p\left(r_\infty-r_L\right)$, we
begin by proving that $r_{\infty} \leq r_L$ for all $L$.  This makes sense at an intuitive level, since $r_L$ is minimizing a function which ``penalizes less" for carrying inventory.  However, in spite of this clear intuition, the analysis is not entirely trivial, as
one function dominating another does not necessarily imply a similar comparison of the appropriate minimizers, and we defer the proof to the technical appendix Section\ \ref{ls-sec-appen}.

\begin{lemma}\label{lslemma-4}
$r_\infty \leq r_L$ for all $L \geq 1$.  
\end{lemma}

Before proceeding, we also derive a certain critical inequality, which we will use to show that the term $h\left(\bbe\left[I^{r_\infty}_L\right]- \bbe\left[I^{r_L}_L\right]\right)$ and the term $p\left(r_\infty - r_L\right)$ essentially ``cancel out".  This inequality follows from the first-order optimality conditions of the convex optimization problem associated with $r_{\infty}$, but requires some care, as the relevant functions are potentially non-differentiable, and the desired statement in principle involves an interchange of expectation and differentiation.  We again defer the proof to the technical appendix Section\ \ref{ls-sec-appen}.
\begin{lemma}\label{lseq-54}
$\sum_{n=1}^\infty \mathbb{P}\left(n r_\infty \geq \sum_{i=1}^n D_{i}\right) \geq \frac{p}{h}.$
\end{lemma}

With Lemmas\ \ref{lslemma-9}, \ref{lslemma-4}, and\ \ref{lseq-54}\ in hand, we now complete the proof of our main result, Theorem\ \ref{lsthm-main}.
\proof{Proof of Theorem \ref{lsthm-main}}
Recall that the remaining term on the right-hand side of (\ref{tobound2}) which we are yet to bound is
\begin{equation}\label{lastdance}
h\left(\bbe\left[I^{r_\infty}_L\right]- \bbe\left[I^{r_L}_L\right]\right) - p\left(r_\infty-r_L\right).
\end{equation}
We first bound $\bbe\left[I^{r_\infty}_L\right]- \bbe\left[I^{r_L}_L\right]$, which by Lemma\ \ref{lslemma-9} equals
$$
\sum_{n=1}^L \frac{1}{n} \bbe\left[\left(nr_\infty - \sum_{i=1}^n D_{i}\right)\mathbb{I}\left(nr_\infty\geq \sum_{i=1}^n D_{i}\right)\right] - \sum_{n=1}^L \frac{1}{n} \bbe\left[\left(nr_L - \sum_{i=1}^n D_{i}\right)\mathbb{I}\left(nr_L\geq \sum_{i=1}^n D_{i}\right)\right].
$$
Combining with Lemma\ \ref{lslemma-4} (i.e. the fact that $r_L \geq r_{\infty}$), which implies that
$$\left(nr_L - \sum_{i=1}^n D_{i}\right)\mathbb{I}\left(nr_L \geq \sum_{i=1}^n D_{i}\right)\geq\ \left(nr_L - \sum_{i=1}^n D_{i}\right)\mathbb{I}\left(nr_\infty \geq \sum_{i=1}^n D_{i}\right),
$$
we conclude that $\bbe\left[I^{r_\infty}_L\right]- \bbe\left[I^{r_L}_L\right]$
is at most
$$
\sum_{n=1}^L \frac{1}{n} \bbe\left[\left(nr_\infty - \sum_{i=1}^n D_{i}\right)\mathbb{I}\left(nr_\infty\geq \sum_{i=1}^n D_{i}\right)\right] -  \sum_{n=1}^L \frac{1}{n} \bbe\left[\left(nr_L - \sum_{i=1}^n D_{i}\right)\mathbb{I}\left(nr_\infty\geq \sum_{i=1}^n D_{i}\right)\right],$$
which itself equals $-\left(r_L - r_\infty \right)\sum_{n=1}^L  \mathbb{P}\left(nr_\infty\geq \sum_{i=1}^n D_{i}\right).$  It follows that (\ref{lastdance}) is at most
\begin{equation}\label{finaleq1}
\left(r_L-r_\infty \right) \left(p - h \sum_{n=1}^L  \mathbb{P}\left(nr_\infty\geq \sum_{i=1}^n D_{i}\right) \right).
\end{equation}
Note that Lemma\ \ref{lseq-54} implies that
\begin{equation}\label{finaleq2}
\sum_{n=L + 1}^{\infty} \mathbb{P}\left(nr_\infty\geq \sum_{i=1}^n D_{i}\right) \geq \frac{p}{h} -
 \sum_{n=1}^L \mathbb{P}\left(nr_\infty\geq \sum_{i=1}^n D_{i}\right).
\end{equation}
Combining (\ref{finaleq1}) and (\ref{finaleq2}), we conclude that (\ref{lastdance}) is at most $\left(r_L - r_\infty \right) h \sum_{n= L + 1}^{\infty} \mathbb{P}\left(nr_\infty\geq \sum_{i=1}^n D_{i}\right).$  It follows from the well-known Chernoff's inequality (cf. \cite{DS}) that $\mathbb{P}\left(nr_\infty\geq \sum_{i=1}^n D_{i}\right) \leq \gamma^n.$  By summing the associated geometric series, and combining with the fact that by definition $r_L \leq \E[D]$, we conclude that (\ref{lastdance}) is at most $h (\E[D] - r_{\infty})(1 - \gamma)^{-1} \gamma^{L+1}.$  Combining the above bound for (\ref{lastdance}) with Lemma\ \ref{lslemma-9}, plugging into (\ref{tobound2}), and applying the fact that $\textrm{OPT}(L) \geq g$, completes the proof.$\Halmos$
\endproof

\section{Conclusion}\label{ls-sec-conc}
In this paper, we proved that for lost sales models with large lead times, the optimality gap of the simple constant-order policy converges \emph{exponentially fast} to zero as the lead time grows with the other problem parameters held fixed, and derived effective explicit bounds for this optimality gap.  This takes a large step towards answering several open questions of \cite{GKLSS}, who recently proved the asymptotic optimality of the constant-order policy in this setting, but whose bounds on the rate of convergence were impractical.  We also demonstrated that the corresponding rate of exponential decay is at least as fast as the exponential rate of convergence of the expected waiting time in a certain single-server queue to its steady-state value, which we proved to be monotone in the ratio of the lost-sales penalty to the holding cost.  For the special case of exponentially distributed demand, we further computed all expressions appearing in our bound in closed form, and numerically evaluated these bounds, demonstrating good performance for a wide range of parameter values.

This work leaves many interesting directions for future research.  First, it would be interesting to investigate the tightness of our exponential bound, e.g. to determine whether our exponential rate captures the true exponential rate of convergence of the optimality gap of the constant-order policy.  Although one can come up with pathological examples for which this is not true, e.g. discrete demand distributions with probability at least $\frac{p}{p+h}$ at $0$ (for which $Q = 0$, $\pi_0$ is optimal amongst all policies for all $L \geq 0$, yet $\gamma > 0$), we conjecture that under mild assumptions $\gamma$ indeed captures the true rate of convergence of the optimality gap.  Second, it is an open challenge to analyze the performance of more sophisticated policies for lost sales inventory models, e.g. affine policies, order-up-policies, and dual-balancing policies.  Another interesting question involves the formal construction and analysis of ``hybrid" algorithms, which e.g. solve large dynamic programs when $L$ is small and transition to using simpler policies when $L$ is large, or use base-stock policies when $\frac{p}{h}$ is large (relative to $L$) and a constant-order policy when $\frac{p}{h}$ is small.  The analysis of such policies would likely require a more precise understanding of the complexity of computing $r_{\infty}$ to any given accuracy, and a deeper understanding of how to apply the techniques of \cite{Halman12} and \cite{Chen12a} to infinite-horizon problems.  Third, it would be interesting to prove that a similar phenomena occurs (albeit possibly with a different simple policy) in more general inventory settings, e.g. models with non-i.i.d. demand, fixed ordering costs, integrality constraints, multiple suppliers and/or products, and network structure.  On a final note, our results and methodology (combined with that of \cite{GKLSS}) provide a fundamentally new approach to lost sales inventory models with positive lead times.  We believe that our approach, combined with other recent developments in inventory theory (e.g. the efficient solution of related dynamic programs), represents a considerable step towards making these models solvable in practice.  Such progress may ultimately help to free researchers from having to use backlogged demand inventory models as approximations to lost-sales inventory models, even when such an approximation is not appropriate, which has been recognized as a major problem in the inventory theory literature (cf. \cite{BV}).

\section{Appendix}\label{ls-sec-appen}
\subsection{Example: exponentially distributed demand}\label{ls-sec-explicit}
In this section, for the special case of exponentially distributed demand, we further compute all expressions appearing in our bound in closed form, and numerically evaluate them, demonstrating good performance for a wide range of parameter values.  Thus suppose demand is exponentially distributed with rate $\lambda$, i.e. mean $\lambda^{-1}$.  In this case, it is well-known that for $r \in [0,\E[D])$, $\bbe[I_\infty^r]$ is the expected steady-state waiting time in a corresponding $M/D/1$ queue, and equals $\frac{r^2\lambda}{2(1-r\lambda)}$ (cf. \cite{Haigh}).  
For $p,h > 0$, let $\tau_{p,h} \stackrel{\Delta}{=} \sqrt{\frac{h}{2p + h}}$, and $\gamma_{p,h} \stackrel{\Delta}{=} (1 - \tau_{p,h})\exp(\tau_{p,h}),$ where it may be easily demonstrated that $\gamma_{p,h} \in (0,1)$.  In that case, when demand is exponentially distributed, Theorem\ \ref{lsthm-main} is equivalent to the following bound.
\begin{corollary}[Case of exponentially distributed demand]\label{ourmainboundexpo}
Suppose $D$ is exponentially distributed with rate $\lambda$.  Then $C(\pi_{r_{\infty}}) = \lambda^{-1}( \sqrt{h (2p + h)} - h )$, and for all $L \geq 1$,
\begin{equation}\label{ls-eq-example}
\frac{C(\pi_{r_{\infty}})}{\textrm{OPT}(L)} \leq 1 +  \bigg( \tau_{p,h} + (\tau^{-1}_{p,h} - 1)\big( e (L+1) \big)^{-1} \bigg) \big( (1 - \gamma_{p,h}) \log(1 + p h^{-1}) \big)^{-1} \gamma^{L+1}_{p,h}.
\end{equation}
\end{corollary}
\proof{Proof}
Suppose $D$ is exponentially distributed with rate $\lambda$.  It is well-known that in this case, for all $r \in [0,\lambda^{-1})$, $\bbe[I_\infty^r] = \frac{r^2 \lambda}{2 (1 - r \lambda)}$ (cf. \cite{Haigh}).  It follows from (\ref{lseqconst}) that for all $r \in [0,\lambda^{-1})$, $C(\pi_r) = h\frac{r^2\lambda}{2(1 - r\lambda)}+ p \lambda^{-1} - p r,$ and $\frac{d}{dr} C(\pi_r) = \frac{h}{2}\big( (\lambda r - 1)^{-2} - 1 \big) - p$ strictly increases from $-p$ to $\infty$ on $[0,\lambda^{-1})$.  It follows that $r_{\infty}$ must be the unique solution to the equation $\frac{d}{dr} C(\pi_r) = 0$ on $[0,\lambda^{-1})$, and hence by a straightforward calculation that $r_{\infty} = \lambda^{-1} (1 - \tau_{p,h})$, and $C(\pi_{r_{\infty}}) = \lambda^{-1}( \sqrt{h (2p + h)} - h ).$  As $\E[\exp(- \theta D)] = \lambda (\lambda + \theta)^{-1}$ for all $\theta \geq 0$, we conclude that
$\phi(\theta) = \exp\big(\lambda^{-1} (1 - \tau_{p,h}) \theta \big) \lambda (\lambda + \theta)^{-1}.$  It then follows from another straightforward calculation that $\vartheta$ equals the unique solution to $\frac{d}{d\theta} \phi(\theta) = 0$ on $[0,\infty)$, $\vartheta = \tau_{p,h} \lambda (1 - \tau_{p,h})^{-1}$, and $\gamma = \gamma_{p,h} = (1 - \tau_{p,h})\exp(\tau_{p,h}).$  It is easily calculated that $Q= \lambda^{-1} \log(1 + p h^{-1})$, and
$$g = h\int_{0}^Q (Q-x)\lambda\exp\left(-\lambda x\right)dx+ p\int_Q^\infty (x-Q)\lambda\exp\left(-\lambda x\right)dx = 
h \lambda^{-1}\log\left(1 + p h^{-1}\right).$$
Combining the above with another straightforward calculation completes the proof.   $\Halmos$
\endproof

Note that (\ref{ls-eq-example}) does not depend on $\lambda$, which follows from the scaling properties of the exponential distribution.  We now numerically evaluate \eqref{ls-eq-example} under different lost-demand penalty and lead time scenarios, with the holding cost fixed to 1, and present the results in Table \ref{ls-table-ub}.  
\begin{table}[H]
\centering
    \begin{tabular}{|c|c|c|c|c|c|c|c|c|c|}
      \hline
      Evaluation of \eqref{ls-eq-example}   & L=1 & L=4 & L=10 & L=20 & L=30  &  L=50  & L=70  &  L=100 \\
      \hline
      p=1/4  &  2.13  &  1.08  &  1.00  &  1.00  &  1.00   &  1.00  &  1.00   &  1.00  \\
      \hline
      p=1    &  3.36  &  1.89  &  1.15  &  1.01  &  1.00   &  1.00  &  1.00   &  1.00   \\
      \hline
      p=4    &  6.42  &  3.99  &  2.62  &  1.72  &  1.34   &  1.08  &  1.02   &  1.00  \\
      \hline
      p=9    &  12.26 &  6.77  &  4.43  &  3.12  &  2.45   &  1.73  &  1.38   &  1.15   \\
      \hline
      p=39   &  62.26 &  27.60 &  14.86 &  9.62  &  7.62   &  5.75  &  4.75   &  3.81  \\
      \hline
      p=99   &  204.5 &  85.21 &  41.77 &  24.43 &  18.20  &  12.92 &  10.49  &  8.49  \\
      \hline
    \end{tabular}
\caption{When $h = 1$, values of \eqref{ls-eq-example} under different $p$ and $L$.}\label{ls-table-ub}
\end{table}
\vspace{-.3in}
We note that consistent with Lemma\ \ref{lsprop-main1}, when $\frac{p}{h}$ is small (i.e. less than or equal to 1), our bounds demonstrate an excellent performance by the constant-order policy even for lead times as small as 10.  When $\frac{p}{h}$ is moderate (i.e. less than or equal to 9), our bounds demonstrate a similarly good performance for lead times on the order of 70.  Even when $\frac{p}{h}$ is very large, our bounds still imply non-trivial performance guarantees, e.g. the constant-order policy is always within a factor of 4 of optimal when $p = 39$ and $L = 100$.  Note that combining (\ref{ls-eq-example}) with our explicit evaluation of $C(\pi_{r_{\infty}})$ yields tight bounds on $\textrm{OPT}(L)$ whenever (\ref{ls-eq-example}) is close to 1.  For example, for $p = \lambda = 1$, our bounds imply that $\textrm{OPT}(10) \in [.63, .73]$.
\subsection{Proof of Lemma\ \ref{lsprop-main1}}
\proof{Proof of Lemma \ref{lsprop-main1}}
Suppose $\frac{p_1}{h_1}< \frac{p_2}{h_2}$, and let $r^{i}_\infty\in\ \argmin_{0\leq r\leq \bbe[D]}\ \left( h_i\bbe[I^r_\infty]+ p_i\bbe[D]- p_ir \right)$, $i=1,2$.
From the respective optimality of $r^1_{\infty}, r^2_{\infty}$, we conclude that
\begin{eqnarray*}
   & \bbe[I^{r^{1}_\infty}_\infty]+ \frac{p_1}{h_1}\bbe[D]- \frac{p_1}{h_1}r^{1}_\infty\leq\ \bbe[I^{r^{2}_\infty}_\infty]+ \frac{p_1}{h_1}\bbe[D]- \frac{p_1}{h_1}r^{2}_\infty,\\
   & \bbe[I^{r^{2}_\infty}_\infty]+ \frac{p_2}{h_2}\bbe[D]- \frac{p_2}{h_2}r^{2}_\infty\leq\  \bbe[I^{r^{1}_\infty}_\infty]+ \frac{p_2}{h_2}\bbe[D]- \frac{p_2}{h_2}r^{1}_\infty.
\end{eqnarray*}
Summing these two inequalities together implies $\left(\frac{p_2}{h_2}- \frac{p_1}{h_1}\right)\left(r^{2}_\infty-r^{1}_\infty\right)\geq 0.$  It follows that $r^{2}_\infty \geq r^{1}_\infty$, and for all $\theta \geq 0$, 
$\exp(\theta r^2_{\infty}) \bbe[\exp(- \theta D)] \geq \exp(\theta r^1_{\infty}) \bbe[\exp(- \theta D)].$
Combining with the definition of $\gamma$ completes the proof.   $\Halmos$
\endproof
\subsection{Proof of Theorem\ \ref{proxychain1}}
Before providing a formal proof, we first provide an intuitive overview, noting that our proof is very similar to many proofs appearing in the literature in which one proves that a ``small perturbation" of a given MDP has additional nice properties (cf. \cite{Filar}).
We proceed by constructing a sequence of random vectors, one for each sufficiently small $\epsilon > 0$, where the random vector ultimately chosen to satisfy the requirements of the theorem will be an appropriate subsequential weak limit (as $\epsilon \downarrow 0$) of these vectors.  Given such an $\epsilon > 0$, we will pick a sufficiently large time $T_{\epsilon}$ s.t. the expected performance of $\pi^{*,L}$ up to time $T_{\epsilon}$ is ``close" to $\textrm{OPT}(L)$.  We then construct a ``modified Markov chain", which behaves very-much like $\lbrace (\mathbf{x}^{*,L}_t, I^{*,L}_t), t \geq 1 \rbrace$, but has an extra ``time-accounting" dimension (i.e. it will be $(L+2)$-dimensional instead of $(L+1)$-dimensional).  This extra dimension will allow the chain (if initialized in a special restart state) to mimic $\pi^{*,L}$ (started from an empty pipeline and zero inventory) for exactly $T_{\epsilon}$ periods, then mimic a policy which orders nothing until the dimensions corresponding to the pipeline and inventory re-enter the $(\mathbf{0},0)$ state, then enter a dummy state (in which it continues to mimic a policy which orders nothing) and stay for a geometrically distributed number of periods, and finally re-enter the special restart state and repeat the process.  This will yield a regenerative process, which will have finite expected regeneration time ``close" to $T_{\epsilon}$ since the pipeline vector and inventory level are bounded by $M_L$ under $\pi^{*,L}$, and which will be aperiodic due to the geometrically distributed number of visits to the dummy state.  We will thus be able to prove, using standard results from the theory of regenerative processes, that this process has a steady-state distribution.  Furthermore, our construction will ensure that the corresponding random vector has approximately (made precise in terms of $\epsilon$) the desired features laid out in the statement of the theorem.  We will then prove that this sequence of random vectors itself has a subsequential weak limit (along a subsequence of $\epsilon$'s converging to 0), which will satisfy the conditions of the theorem. 
\proof{Proof of Theorem\ \ref{proxychain1}}
Recall that $\textrm{OPT}(L) \geq g > 0$.  Similarly, as $C(\pi_0) = p$, $\textrm{OPT}(L) \leq p$.  Given any $\epsilon \in (0,g)$, we now construct a random vector $(\chi^{\epsilon}_{\infty},\mathcal{I}^{\epsilon}_{\infty})$ 
which satisfies (\ref{enum1}) - (\ref{enum3}) (after replacing $\chi^*$ and ${\mathcal I}^*$ by $\chi^{\epsilon}_{\infty}$ and $\mathcal{I}^{\epsilon}_{\infty}$ in all relevant statements and expressions), and for which
\begin{equation}\label{enum4approx}
\big|h \E[{\mathcal I}^{\epsilon}_{\infty}] +  p \bbe[D_1] - p \bbe[\chi^{\epsilon}_{1,\infty}] - \textrm{OPT}(L)\big| < \epsilon,
\end{equation}
i.e. (\ref{enum4}) is approximately satisfied.  It follows from the definition of $\limsup$ that there exists
$T_{\epsilon} > 2  (L + \frac{M_L}{\E[D]} + 2)\big(h M_L + (\E[D] + 1)p\big) \epsilon^{-1}$ s.t. $(\textrm{OPT}(L) - \frac{\epsilon}{2}) T_{\epsilon} < \sum_{t=1}^{T_{\epsilon}} \bbe\left[C^{\pi^{*,L}}_t \right] < 
(\textrm{OPT}(L) + \frac{\epsilon}{2}) T_{\epsilon}.$  We now construct an $(L+2)$-dimensional non-negative discrete-time Markov process $\lbrace \mathbf{Y}^{\epsilon}_t, t \geq 1 \rbrace = \lbrace (\chi^{\epsilon}_t, {\mathcal I}^{\epsilon}_t, \tau^{\epsilon}_t) , t \geq 1 \rbrace$, where $\chi^{\epsilon}_t = (\chi^{\epsilon}_{1,t},\ldots,\chi^{\epsilon}_{L,t})$ is an $L$-dimensional random vector, and ${\mathcal I}^{\epsilon}_t$ and $\tau^{\epsilon}_t$ are both 1-dimensional r.v.s.  Let $\lbrace B_t, t \geq 1 \rbrace$ denote an i.i.d. sequence of Bernoulli r.v.s, each of which equals 1 w.p. $\frac{1}{2}$ and 0 w.p. $\frac{1}{2}$.  Then $\lbrace \mathbf{Y}^{\epsilon}_t, t \geq 1 \rbrace$ evolves as follows.  $\chi^{\epsilon}_1 = \mathbf{0}, {\mathcal I}^{\epsilon}_1 = 0, \tau^{\epsilon}_1 = 1$.  For $t \geq 1$, the chain evolves as follows.  $\chi^{\epsilon}_{i,t+1} = \chi^{\epsilon}_{i+1,t}$ for $i \in [1,L-1]$, and ${\mathcal I}^{\epsilon}_{t+1} = ({\mathcal I}^{\epsilon}_t + \chi^{\epsilon}_{1,t} - D_t)^+$.  If $\tau^{\epsilon}_t \in [1,T_{\epsilon})$, then $\chi^{\epsilon}_{L,t+1} = f^{\pi^{*,L}}(\chi^{\epsilon}_t,{\mathcal I}^{\epsilon}_t)$ and $\tau^{\epsilon}_{t+1} = \tau^{\epsilon}_t + 1$.  If $\tau^{\epsilon}_t  = T_{\epsilon}$ and $(\chi^{\epsilon}_t,{\mathcal I}^{\epsilon}_t) \neq (\mathbf{0},0)$, then $\chi^{\epsilon}_{L,t+1} = 0$ and $\tau^{\epsilon}_{t+1} = T_{\epsilon}$.  If $\tau^{\epsilon}_t  = T_{\epsilon}$ and $(\chi^{\epsilon}_t,{\mathcal I}^{\epsilon}_t) = (\mathbf{0},0)$, then $\chi^{\epsilon}_{L,t+1} = 0$ and $\tau^{\epsilon}_{t+1} = 0.$  If $\tau^{\epsilon}_t = 0$, then $\chi^{\epsilon}_{L,t+1} = 0$ and $\tau^{\epsilon}_{t+1} = B_t$.  One may easily verify the following properties of $\lbrace \mathbf{Y}^{\epsilon}_t, t \geq 1 \rbrace$.  Let $D'_1$ denote another r.v. distributed as ${\mathcal D}$ (independent of $\lbrace D_t, t \geq 1 \rbrace$), $z(x) \stackrel{\Delta}{=} \E[ h(x - D'_1)^+ + p(x - D'_1)^-]$, and $\hat{T}_{\epsilon}$ denote a r.v. distributed as the time between the chain's initial and second visit to $(\mathbf{0},0,1)$.
\begin{enumerate}[(a)]
\item \label{firsta} It follows from a straightforward induction (essentially Lindley's recursion, cf. \cite{GKLSS}) that for all $t \geq 1$, $\mathcal{I}^{\epsilon}_{t+L} \sim \max_{j=0,\ldots,L}\left(\sum_{i=1}^j(\chi^{\epsilon}_{L+1-i,t} - D_{t + L- i}) + \delta_{j,L} \mathcal{I}^{\epsilon}_t \right)$.
\item \label{firstb} It follows directly from the Markov chain dynamics that for all $t \geq 1$, $\E[\mathcal{I}^{\epsilon}_{t+1}] = \E[({\mathcal I}^{\epsilon}_t + \chi^{\epsilon}_{1,t} - D_t)^+]$, and thus $\E[({\mathcal I}^{\epsilon}_t + \chi^{\epsilon}_{1,t} - D_t)^-] = 
\E[\mathcal{I}^{\epsilon}_{t+1}] - \E[\mathcal{I}^{\epsilon}_t] + \E[D_t] - \E[\chi^{\epsilon}_{1,t}]$;
\item \label{firstc} It follows directly from the Markov chain dynamics that for all $t \geq 1$,
$\chi^{\epsilon}_{i,t+1} \sim \chi^{\epsilon}_{i+1,t}$ for $i \in [1,L-1]$.
\item \label{firstd} Conditional on the event $\lbrace \mathbf{Y}^{\epsilon}_t = (\mathbf{0},0,1) \rbrace$, the joint distribution of $\lbrace \mathbf{Y}^{\epsilon}_i , i \in [t, t  + T_{\epsilon}-1] \rbrace$ is identical to that of $\lbrace (\mathbf{x}^{*,L}_i, I^{*,L}_i,i), i \in [1,T_{\epsilon}] \rbrace$.
\item \label{firste} Conditional on the event $\lbrace \tau^{\epsilon}_t = T_{\epsilon} \rbrace$, the expected number of time steps until the Markov chain next enters the state $(\mathbf{0},0,0)$ is at most $L + \frac{M_L}{\E[D]}$.  We note that this follows from (\ref{isbounded}), combined with a straightforward application of Wald's identity and the fact that ordering nothing for $L$ periods clears the pipeline vector.
\item \label{firstf} Upon entering state $(\mathbf{0},0,0)$, the Markov chain remains in that state for a geometrically distributed number of time steps (with mean 2), and then transitions to state $(\mathbf{0},0,1)$.
\item \label{boundTeps} W.p.1, $\hat{T}_{\epsilon} \geq T_{\epsilon}$, and $\E[\hat{T}_{\epsilon}] - T_{\epsilon} \leq L + \frac{M_L}{\E[D]} + 2$.
\item \label{zboundeps} 
$
0 \leq 
\E[\sum_{t=1}^{\hat{T}_{\epsilon}} z(\chi^{\epsilon}_{1,t} + \mathcal{I}^{\epsilon}_t)]
-
\E[\sum_{t=1}^{T_{\epsilon}} z(\chi^{\epsilon}_{1,t} + \mathcal{I}^{\epsilon}_t)]
\leq 
(L + \frac{M_L}{\E[D]} + 2)(h M_L + p \E[D])$.
\end{enumerate}
Combining (\ref{firstd}) - (\ref{firstf}) with the basic definitions associated with the theory of regenerative processes (here we refer the interested reader to \cite{Asmussen} for an excellent overview), we conclude that 
$\lbrace \mathbf{Y}^{\epsilon}_t, t \geq 1 \rbrace$ is a discrete-time aperiodic regenerative process, with regeneration points coinciding with visits to $(\mathbf{0},0,1)$.  It thus follows from standard results in the theory of regenerative processes (cf. \cite{Asmussen}) that $\lbrace \mathbf{Y}^{\epsilon}_t, t \geq 1 \rbrace$ converges weakly (as $t \rightarrow \infty$) to a limiting random vector $\mathbf{Y}^{\epsilon}_{\infty} = (\chi^{\epsilon}_{\infty}, \mathcal{I}^{\epsilon}_{\infty},\tau^{\epsilon}_{\infty})$.  That $(\chi^{\epsilon}_{\infty}, \mathcal{I}^{\epsilon}_{\infty})$ satisfies 
(\ref{enum1}) - (\ref{enum3}) then follows from (\ref{firsta}) - (\ref{firstc}).  We now prove that $(\chi^{\epsilon}_{\infty}, \mathcal{I}^{\epsilon}_{\infty})$ also satisfies (\ref{enum4approx}).  Again applying the standard theory of regenerative processes (cf. \cite{Asmussen}), we conclude that $\E[z(\chi^{\epsilon}_{1,\infty} + \mathcal{I}^{\epsilon}_{\infty})] = \frac{\E[\sum_{t=1}^{\hat{T}_{\epsilon}} z(\chi^{\epsilon}_{1,t} + \mathcal{I}^{\epsilon}_t)]}{\E[\hat{T}_{\epsilon}]}$.  Combining with (\ref{firsta}) - (\ref{zboundeps}), the definition of $T_{\epsilon}$, and some straightforward algebra completes the proof of (\ref{enum4approx}), and we omit the details.   To complete the proof of the theorem, we observe that the sequence of probability measures corresponding to $\lbrace (\chi^{\frac{1}{n}}_{\infty}, \mathcal{I}^{\frac{1}{n}}_{\infty}), n \geq \lceil g^{-1} \rceil \rbrace$ is tight since all associated random vectors have all components bounded in absolute value w.p.1 by $M_L$, where we refer the interested reader to \cite{Bill13} for a review of tightness and related notions.  Thus this sequence of measures has at least one subsequence which converges weakly to some subsequential weak limit $(\chi^*,{\mathcal I}^*)$, which by the already proven properties of $(\chi^{\frac{1}{n}}_{\infty}, \mathcal{I}^{\frac{1}{n}}_{\infty})$ and definition of weak limit will satisfy all required conditions of the theorem.  $\Halmos$.
\endproof

\subsection{Proof of Lemma\ \ref{lslemma-4}}
\proof{Proof of Lemma\ \ref{lslemma-4}}
Suppose for contradiction that there exists $L \in [1,\infty)$ s.t. $r_L< r_\infty$.  Note that in this case, both $r_L, r_{\infty} < \E[D]$, and thus by Lemma\ \ref{lslemma-9} both $\E[I^{r_{\infty}}_{\infty}], \E[I^{r_L}_{\infty}] < \infty$.  From definitions and the associated respective optimality of $r_L,r_{\infty}$, we conclude that
\begin{eqnarray*}
   & h\bbe[I^{r_\infty}_\infty]+ p\bbe[D]- pr_\infty\leq\  h\bbe[I^{r_L}_\infty]+ p\bbe[D]- pr_L,\\
   & h\bbe[I^{r_L}_L]+ p \bbe[D]- p r_L\leq\  h\bbe[I^{r_\infty}_L]+ p\bbe[D]- p r_\infty.
\end{eqnarray*}
Summing these two inequalities together implies that
$$
   \bbe[I^{r_\infty}_\infty]+ \bbe[I^{r_L}_L]\leq\ \bbe[I^{r_L}_\infty]+ \bbe[I^{r_\infty}_L],
$$
which, by Lemma\ \ref{lslemma-9}, is equivalent to
\begin{equation}\label{lseq-a1}
   \sum_{n=L+1}^\infty \frac{1}{n} \bbe\left[\left(nr_\infty - \sum_{i=1}^n D_{i}\right)^+\right]\leq\ \sum_{n=L+1}^\infty \frac{1}{n} \bbe\left[\left(nr_L - \sum_{i=1}^n D_{i}\right)^+\right] < \infty.
\end{equation}
Here we note that the desired result intuitively follows from (\ref{lseq-a1}) and the monotonicity of the relevant functions, i.e. the fact that $x > y$ implies $\bbe\left[\left(n x - \sum_{i=1}^n D_{i}\right)^+\right] \geq \bbe\left[\left(n y - \sum_{i=1}^n D_{i}\right)^+\right]$.  However, we must rule out certain subtle problems that could potentially arise from the function $\bbe\left[\left(n r - \sum_{i=1}^n D_{i}\right)^+\right]$ not being \emph{strictly} monotonic in $r$, and proceed as follows.  
Definitions, non-negativity, and the fact that $r_L < r_{\infty}$, together imply that
\begin{eqnarray*}
\left(nr_\infty - \sum_{i=1}^n D_{i}\right)^+ &=& \mathbb{I}\left(nr_{\infty} \geq \sum_{i=1}^n D_{i}\right)\left(nr_{\infty} - \sum_{i=1}^n D_{i}\right)
\\&\geq& \mathbb{I}\left(n r_L \geq \sum_{i=1}^n D_{i}\right)\left(nr_{\infty} - \sum_{i=1}^n D_{i}\right).
\end{eqnarray*}
Combining with (\ref{lseq-a1}), we conclude that
$$\sum_{n=L+1}^\infty \frac{1}{n} \bbe\left[ \mathbb{I}\left(n r_L \geq \sum_{i=1}^n D_{i}\right)\left(nr_{\infty} - \sum_{i=1}^n D_{i}\right) \right] \leq \sum_{n=L+1}^\infty \frac{1}{n} \bbe\left[ \mathbb{I}\left(n r_L \geq \sum_{i=1}^n D_{i}\right)\left(n r_L - \sum_{i=1}^n D_{i}\right) \right] < \infty.$$
It follows that
$$\sum_{n=L+1}^\infty \frac{1}{n} \bbe\left[ \mathbb{I}\left(n r_L \geq \sum_{i=1}^n D_{i}\right)\left( n r_{\infty} - n r_L \right) \right] \leq 0.$$
However, since by assumption $n r_{\infty} - n r_L > 0$, it follows from non-negativity that
$$\sum_{n=L+1}^\infty \bbe\left[ \mathbb{I}\left(n r_L \geq \sum_{i=1}^n D_{i}\right) \right] = 0,$$
and thus $\p\left(n r_L \geq \sum_{i=1}^n D_{i}\right) = 0$ for all $n \geq L + 1$.  Further noting that
$\mathbb{P}\left(nr_L\geq \sum_{i=1}^n D_{i}\right)\geq\ \mathbb{P}^n\left(r_L\geq D_1\right)$, we conclude that $\mathbb{P}\left(r_L\geq D_1\right)= 0$.  It follows that $\bbe[I^{r_L}_\infty]= \bbe[I^{r_L}_L]= 0$, and thus by (\ref{lseqconst}), $C\left(\pi_{r_L}\right) = C_L\left(r_L\right)$.  Combining with Lemma\ \ref{lslemma-0}, which implies that $C_L\left(r_L\right) \leq C\left(\pi_{r_\infty}\right)$, and the optimality of $r_{\infty}$, we conclude that $r_L \in \argmin_{0\leq r\leq \bbe[D]} C\left(\pi_{r}\right)$.  However, as $r_{\infty}$ is by definition the infimum of $\argmin_{0\leq r\leq \bbe[D]} C\left(\pi_{r}\right)$, the fact that $r_L < r_{\infty}$ thus yields a contradiction, completing the proof.  $\Halmos$
\endproof
\subsection{Proof of Lemma\ \ref{lseq-54}}
\proof{Proof of Lemma\ \ref{lseq-54}}
Since $r_\infty<\bbe[D]$, there exists $\delta > 0$ s.t. $r_\infty+\epsilon < \bbe[D]$ for all $\epsilon \in [0,\delta]$.  Let us fix any such $\epsilon > 0$.  The definition and associated optimality of $r_\infty$ implies that $C(\pi_{r_{\infty}}) \leq C(\pi_{r_{\infty} + \epsilon})$.  Combining with Lemma\ \ref{lslemma-9} and (\ref{lseqconst}), we conclude that
$$h\sum_{n=1}^\infty \frac{1}{n} \bbe\left[\mathbb{I}\left(nr_\infty\geq \sum_{i=1}^n D_{i}\right)\left(nr_\infty - \sum_{i=1}^n D_{i}\right)\right]$$
is at most
$$h\sum_{n=1}^\infty \frac{1}{n} \bbe\left[\mathbb{I}\left(n(r_\infty+\epsilon)\geq \sum_{i=1}^n D_{i}\right)\left(n(r_\infty+\epsilon)- \sum_{i=1}^n D_{i}\right)\right] - p \epsilon.
$$
Combining with the fact that
$$
   \mathbb{I}\left(nr_\infty\geq \sum_{i=1}^n D_{i}\right)\left(nr_\infty- \sum_{i=1}^n D_{i}\right) \geq  \mathbb{I}\left(n(r_\infty+\epsilon)\geq \sum_{i=1}^n D_{i}\right)\left(nr_\infty  - \sum_{i=1}^n D_{i}\right),
$$
it follows that
$$h\sum_{n=1}^\infty \frac{1}{n} \bbe\left[\mathbb{I}\left(n(r_\infty+\epsilon)\geq \sum_{i=1}^n D_{i}\right)\left(nr_\infty  - \sum_{i=1}^n D_{i}\right)\right]$$
is at most
$$h\sum_{n=1}^\infty \frac{1}{n} \bbe\left[\mathbb{I}\left(n(r_\infty+\epsilon)\geq \sum_{i=1}^n D_{i}\right)\left(n(r_\infty+\epsilon)- \sum_{i=1}^n D_{i}\right)\right] - p \epsilon.
$$
Equivalently (as all relevant sums are finite)
$$\sum_{n=1}^\infty \bbp\left(n(r_\infty+\epsilon)\geq \sum_{i=1}^n D_{i}\right) \geq \frac{p}{h}.
$$
As this holds for all sufficiently small $\epsilon$, the only remaining step is to demonstrate validity at $\epsilon = 0$.
By monotonicity, for each fixed $n$ and all $\epsilon \in [0,\delta]$,
$$\bbp\left(n(r_\infty+\epsilon)\geq \sum_{i=1}^n D_{i}\right)
\leq
\bbp\left(n(r_\infty + \delta)\geq \sum_{i=1}^n D_{i}\right).
$$
Furthermore, since $r_{\infty} + \delta < \E[D]$, for any fixed $\nu > 0$, there exists $M_{\nu} < \infty$ (depending only on $\nu, {\mathcal D}, r_{\infty}, \delta$) s.t.
$$\sum_{n=M_{\nu}+1}^\infty \bbp\left(n(r_\infty+\delta) \geq \sum_{i=1}^n D_{i}\right) \leq \nu.
$$
Indeed, the above follows from a standard argument (the details of which we omit) in which each term is bounded using Chernoff's inequality, and the terms are summed as an infinite series (cf. \cite{DS}).  Combining the above, we conclude that for all $\nu > 0$, and $\epsilon \in (0,\delta]$,
$$\sum_{n=1}^{M_{\nu}} \bbp\left(n(r_\infty+\epsilon)\geq \sum_{i=1}^n D_{i}\right) \geq \frac{p}{h} - \nu.
$$
As $\bbp\left(n x \geq \sum_{i=1}^n D_{i}\right)$ is a right-continuous function of $x$  (by the right-continuity of cumulative distribution functions), it follows that $\sum_{n=1}^{M_{\nu}} \bbp\left(n x \geq \sum_{i=1}^n D_{i}\right)$ is similarly right-continuous in $x$.  Right-continuity at $\epsilon = 0$ follows, and we conclude that
$$\sum_{n=1}^{M_{\nu}} \bbp\left(n r_\infty \geq \sum_{i=1}^n D_{i}\right) \geq \frac{p}{h} - \nu.
$$
As this holds for all $\nu$, letting $\nu \downarrow 0$ completes the proof.$\Halmos$
\endproof

\section*{Acknowledgments}
The authors thank Maury Bramson, Jim Dai, Bruce Hajek, Nir Halman, Tim Huh, Ganesh Janakiraman, Dmitriy A. Katz-Rogozhnikov, Retsef Levi, Yingdong Lu, Marty Reiman, Alan Scheller-Wolf, Mayank Sharma, Mark Squillante, Kai Wang, Jiawei Zhang, and Paul Zipkin for several stimulating discussions.  David A. Goldberg also gratefully acknowledges support from NSF grant no. 1453929.

\bibliographystyle{nonumber}

\end{document}